\documentclass[runningheads]{llncs}
\usepackage{graphicx}
\usepackage{color}
\usepackage{amsmath}
\usepackage{amssymb}
\usepackage{cite}
\usepackage{hyperref}
\usepackage{algorithm}
\usepackage{algorithmic}
\usepackage{todonotes}

\newcommand{\ZZ}{\mathbb{Z}} 
\newcommand{\RR}{\mathbb{R}} 

\begin{document}
\title{Online Learning for Scheduling MIP Heuristics
\thanks{This work was partially funded by the Deutsche Forschungsgemeinschaft (DFG, German Research Foundation) under Germany’s Excellence Strategy – The Berlin Mathematics Research Center MATH+ (EXC-2046/1, 390685689) and by the German Federal Ministry of Education and Research (BMBF) within the Research Campus MODAL (05M14ZAM, 05M20ZBM) and supported by a Google Research Award.}}
%
%
\author{Antonia Chmiela\inst{1},
Ambros Gleixner\inst{1,2},
Pawel Lichocki\inst{3},
Sebastian Pokutta\inst{1,4}}
%
\authorrunning{A. Chmiela et al.}
%
\institute{Zuse Institute Berlin, 
\email{\{chmiela,gleixner,pokutta\}@zib.de} \and
Hochschule f{\"u}r Technik und Wirtschaft Berlin \and
Google Research, \email{pawell@google.com} \and
Technische Universit{\"a}t Berlin}
\maketitle              
\begin{abstract}
Mixed Integer Programming (MIP) is NP-hard, and yet modern solvers often solve 
large real-world problems within minutes.
This success can partially be attributed to heuristics.
Since their behavior is highly instance-dependent, relying on hard-coded rules derived from empirical testing on a large heterogeneous corpora of benchmark instances might lead to sub-optimal performance.
In this work, we propose an online learning approach that adapts the application of heuristics towards the single instance at hand.
We replace the commonly used static heuristic handling with an adaptive framework exploiting past observations about the heuristic's behavior to make future decisions.
In particular, we model the problem of controlling Large Neighborhood Search and Diving -- two broad and complex classes of heuristics -- as a multi-armed bandit problem.
Going beyond existing work in the literature, we control two different classes of heuristics simultaneously by a single learning agent.
We verify our approach numerically and show consistent node reductions over the MIPLIB 2017 Benchmark set.
For harder instances that take at least 1000 seconds to solve, we observe a speedup of $4\%$.

\keywords{Mixed Integer Programming \and Machine Learning \and Heuristics.}
\end{abstract}
\section{Introduction}
\label{sec:introduction}

A multitude of problems arising from real-world applications can be modeled as \emph{Mixed Integer Problems (MIPs)}.
Because of that, there is high interest in finding ways to solve MIPs efficiently.
Generally, the \emph{Branch-and-Bound (B\&B)} framework \cite{land60} is used which decomposes the optimization problem in smaller subproblems that are then easier to handle.
Since this approach involves a variety of decisions that significantly influence its behavior, the idea of using machine learning (ML) has gained interest:
ML has been used to find good solver parameters \cite{hutter09, hutter11, dambrosio20}, to improve node \cite{he14, yilmaz21}, variable \cite{khalil16, lodi17, balcan2018learning, etheve20, nair20, scavuzzo22}, and cut selection \cite{baltean19, tang20, huang22, turner22, paulus22}, and to detect decomposable structures \cite{kruber17}.

The objective of B\&B is to solve MIPs to global optimality.
However,
it is often not feasible to wait until the optimum is found, thus finding good 
feasible solutions early on is important.
Primal heuristics are crucial for this: 
In \cite{berthold132}, the authors showed that heuristics improved the primal bound by $80\%$ and the solving time by $30\%$ on average.
An overview of different primal heuristics and their impact can be found in \cite{lodi17, berthold13, berthold18}.

Primal heuristics are powerful but can be very costly, thus it is important to be strategic about how they are applied in practice.
Controlling their behavior by hard-coded rules derived from empirical tests on heterogeneous benchmark sets leads to strategies that work averagely well on a broad variety of instances.
However, since the performance of heuristics is highly instance-dependent, this might lead to suboptimal behavior.
For example, primal performance can be significantly improved by deriving problem-specific heuristic settings \cite{chmiela21}. 

In this work, we present an online learning approach to control primal heuristics within B\&B. 
We model heuristic selection as a multi-armed bandit problem and exploit past observations of heuristics' behavior to learn on-the-fly which heuristics are most likely to be successful.
Our scheduler is, thus, capable to adapt to and leverage specific characteristic of the problem at hand.
In particular, we control Large Neighborhood Search and Diving, two significantly different and complex classes of heuristics. \\

\noindent\textbf{Contribution.}
To the best of our knowledge, this is the first time when two different classes of heuristics are treated simultaneously by a single learning agent.
To summarize:
\begin{enumerate}
    \item We propose an \textbf{online learning approach for heuristic 
    scheduling} to replace more static heuristic handling 
    (Section~\ref{sec:scheduler}),  
    \item We support our findings by \textbf{extensive computational 
    experiments} on a heterogeneous benchmark test set to numerically verify 
    our approach (Section~\ref{sec:results}).
\end{enumerate}

\paragraph{Related Work.}
Since heuristics have a large impact on solver performance, using ML to develop 
new strategies and to optimize their usage is a topic of ongoing research.
For instance, \cite{nair20} use neural networks to derive variable assignments to find primal solutions.
The authors in \cite{chen22} propose a bi-layer prediction model utilizing graph convolutional networks designed to help heuristics find solutions faster.
To improve the usage of heuristics, the authors in \cite{khalil16} learn an oracle that aims to predict at which nodes a heuristic will be successful or not.
Whereas in \cite{chmiela21} a data-driven heuristic scheduling framework is proposed that learns problem-specific heuristic schedules to find many solutions at minimal cost.

In \cite{hendel18, hendel182} adaptive heuristics are built that use bandit algorithms to decide which heuristics to additionally run.
In particular, their ALNS heuristic \cite{hendel18} inspired the framework we present here: 
While ALNS was designed as another primal heuristic to be added in the pool of available heuristics, we extend it to a framework that aims to replace static heuristic handling and that can be easily extendable to handle any class of heuristics.

\section{Background}
\label{sec:mips}

\paragraph{Mixed Integer Problems.}
A MIP is an optimization problem of the form
\begin{align} \label{mip} \tag{P}
\min_x \ c^T x, \text{ s.t. }  Ax \leq b, \ x_i \in \ZZ, i \in I,
\end{align}
with matrix $A \in \RR^{m \times n}$, vectors $b,c \in \RR^m$ and index set $I 
\subseteq [n]$.
To solve \eqref{mip}, B\&B partitions the feasible region, resulting in a tree structure with nodes correspond to the simpler subproblems.

\emph{Primal Heuristics.}
Heuristics aim to find feasible solutions for \eqref{mip}.
Generally, a solver utilizes a variety of heuristics exploiting different ideas to find high-quality solutions.
Two of the most complex and time consuming classes of heuristics are \emph{diving} and \emph{Large Neighborhood Search (LNS)}.
Diving heuristics examine a single probing path by sub-sequentially fixing variables according to a specific rule.
In contrast, LNS builds a neighborhood around a reference point by fixing a certain percentage of variables and then solving the resulting sub-MIP.
Since no heuristic is guaranteed to be successful, the solver iterates over all available heuristics in a predefined order to hopefully find a new solution.
Good heuristics, like diving and LNS, are typically computationally expensive. 
Thus, it is especially important to be strategic about controlling them.

\emph{Multi-Armed Bandit Problem.}
Given a set of actions $\mathcal{A}$, an agent aims to select a series of actions with maximal cumulative reward.
In every iteration $t$, an action $a_t \in \mathcal{A}$ is selected for which a reward $r_t \in [0,1]$ is observed.
Since the agent only learns how the selected action behaves, a good strategy entails a balance between exploring unknown actions and exploiting the ones that performed well in the past.
There are various approaches to finding a good strategy, see \cite{bubeck12, lattimore20}.

\section{Scheduling Primal Heuristics Online}
\label{sec:scheduler}

We present an online learning approach that models heuristic handling as a multi-armed bandit problem.
Thereby, the set of actions $\mathcal{A}$ corresponds to the set of heuristics $\mathcal{H}$ we want to control.
Two main challenges arise when modeling the scheduling of heuristics this way:
(i) defining a suitable reward function and (ii) choosing the right bandit algorithm.
After presenting our online scheduling framework, we describe how we tackle 
both in Section~\ref{sec:reward} and \ref{sec:bandit}, respectively.

\subsection{The Online Scheduling Framework}

\begin{figure}[t]
\centering
\includegraphics[width=\textwidth]{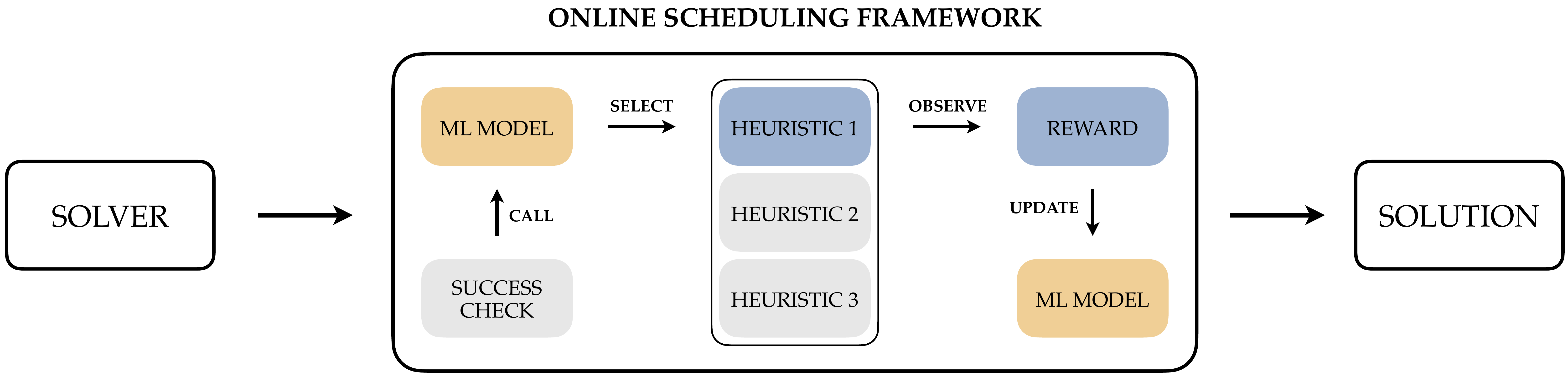}
\caption{\textbf{Visualisation of the Online Scheduling Framework:} When the solver decides to run heuristics, it is checked if the scheduler was successful enough in the past. If so, a bandit algorithm selects a heuristic which is executed with specific working limits. A reward is observed and then used to update the bandit as well as the working limits. A solution is returned to the solver if one was found.}
\label{fig:scheduler}
\end{figure}

The scheduler controls a set of primal heuristics $\mathcal{H}$.
Each heuristic has special working limits influencing its behavior.
Whenever the scheduler is called, it selects and executes one heuristic $h \in \mathcal{H}$.
Depending on how $h$ performed, we dynamically adapt some of its working limits.
This way, we not only tailor heuristic handling to the instance at hand but also reduce the number of user-defined parameters.
To summarize, the scheduler executes the following steps:
\begin{enumerate}
    \item \textbf{Select heuristic} $h$ using a suitable bandit algorithm (introduced in Section~\ref{sec:bandit}).
    \item \textbf{Execute heuristic} $h$ using the current working limits.
    \item \textbf{Observe reward} $r$ after executing $h$ (introduced in Section~\ref{sec:reward}).
    \item \textbf{Update bandit algorithm and working limits} of $h$ using reward $r$. 
\end{enumerate}
An overview of the scheduling framework is shown in Figure~\ref{fig:scheduler}.

Often, a solver finds an optimal solution noticeably faster than it proves the solution's optimality \cite{berthold17}.
Thus, always running heuristics with the same frequency is not necessarily the best strategy.
To dynamically adapt how often the scheduler is executed, we track how often no solution was found.
Whenever it is unsuccessful for too long, we skip a number of future calls to the scheduler:
We skip $\lfloor \exp(\beta n_{\text{fail}}) \rfloor - 1$ calls, where $n_{\text{fail}}$ counts consecutively failed calls and $\beta = 0.1$.

At the beginning of the solving process, when heuristics run for the first time, the scheduler does not have any information about the heuristic's behavior yet.
Thus, any bandit algorithm would start by selecting heuristics at random.
To avoid uninformed decisions, our framework uses expert knowledge to warmstart the bandit strategy:
We execute all heuristics in their default order first and observe their rewards;
only then the bandit algorithm takes over.

This is a general heuristic scheduling framework that can be applied to an arbitrary set $\mathcal{H}$.
However, as mentioned in Section~\ref{sec:mips}, we focus on LNS and diving since they cover the majority of the more complex heuristics.
We control different types of working limits directly influencing the cost and success probability of the heuristics:
For LNS, we impose a target fixing rate and for diving, we control the LP 
resolve frequency.
We adapt both as follows.

The target fixing rate controls how many integer variable should be fixed in 
the sub-MIP.
This directly influences the success rate as well as the costs of the heuristic.
The more variables are fixed, the easier but also the more restrictive the resulting subproblem becomes.
To dynamically adapt the fixing rate, we use the same approach as presented in \cite{hendel18}.
Let us denote by $f_h^t \in [0,1]$ the target fixing rate of heuristic $h$ at iteration $t$.
Assuming that $h$ was selected, we have
\[
f_h^{t+1} =
\begin{cases}
    \max\{(1 - \gamma) f_h^t, f_{min} \}, & \text{if } h \text{ found solution} \text{ or sub-MIP was infeasible,} \\
    \min\{(1 + \gamma) f_h^t, f_{max} \}, & \text{otherwise,}
\end{cases}
\]
with factor $\gamma \in [0,1]$ and target fixing rate limits $f_{min}, f_{max} \in [0,1]$.
We choose $\gamma = 0.1$, $f_{min} = 0.3$, $f_{max} = 0.9$, and $f^0_h = f_{max}$ for all LNS heuristics.

Diving heuristics successively fix integer variables and reoptimize the corresponding LP relaxation in between.
If the LP is solved more often, diving becomes more expensive, but also more successful:
Fixings that led to infeasibility can be detected earlier and then be corrected by backtracking.
To control how often the LP is solved, the fraction of variables is tracked 
that had their domains changed since the last LP solve.
If this fraction exceeds a threshold parameter $q$, an LP solve is triggered; larger $q$ results in less frequent LP solves. 

We dynamically adjust this threshold in a similar fashion to the target fixing rate of LNS heuristics.
Let us denote by $q_h^t \in [0,1]$ the value for diving heuristic $h$ at iteration $t$.
If $h$ was selected at $t$, then
\[
q_h^{t+1} =
\begin{cases}
    \max\{(1 - \eta) q_h^t, q_{min} \}, & \text{if } h \text{ did not find an incumbent at iteration } t \\
    \min\{(1 + \eta) q_h^t, q_{max} \}, & \text{otherwise}
\end{cases}
\]
for factor $\eta \in [0,1]$ and the bounds $q_{min}, q_{max} \in [0,1]$.
Thus, if $h$ was not successful, we reduce $q_h^t$ to increase the success probability of $h$ in the future.
Otherwise, we increase the value to reduce the cost of executing $h$.
We choose $\eta = 0.1$, $q_{min} = 0.05$, $q_{max} = 0.3$, and $q_h^0 = 
q_{min}$ for all diving heuristics.

\subsection{Choosing a Reward Function}
\label{sec:reward}

The simplest choice to reward a heuristic $h \in \mathcal{H}$ would be the binary function
\[
r_{\text{sol}}(h, t) =
\begin{cases}
    1, & \text{if } h \text{ found an incumbent at iteration } t \\
    0, & \text{otherwise}.
\end{cases}
\]
However, heuristics find improving solutions rather rarely: For instance, on the test set we consider in our experiments, the default settings of the solver found on average only 12 incumbents.
Thus, using $r_{\text{sol}}$ as the only reward signal might not give enough feedback to the agent.
Furthermore, $r_{\text{sol}}$ lacks a lot of important information.
For example, a heuristic that fails fast is preferable over one that needs more time to terminate without a solution.
Furthermore, if a solution is found, its quality should also be considered.
Besides the obvious preference for better solutions, considering the current stage of the solving process is vital: 
At the beginning, it is much easier to find a new incumbent than at a more advanced stage.
Another problem is that $r_{\text{sol}}$ implicitly assumes the only objective 
of heuristics is finding solutions.
This is not always true, for instance, diving heuristics can also generate conflict constraints \cite{achterberg07}, which profits the solver.

Thus, besides $r_{\text{sol}}$, we consider three additional metrics to reward $h$:
\begin{enumerate}
    \item $r_{\text{gap}}$ to reward the quality of the new incumbent if $h$ was successful,
    \item $r_{\text{eff}}$ to punish the effort it took to execute $h$,
    \item $r_{\text{conf}}$ to reward the number of conflict constraints $h$ found.
\end{enumerate}
The overall reward function $r$ is then
\[
r(h,t) = \lambda_1 r_{\text{sol}}(h,t) + \lambda_2 r_{\text{gap}}(h,t) + \lambda_3 r_{\text{eff}}(h,t) + \lambda_4 r_{\text{conf}}(h,t),
\]
with $\lambda_i \in [0,1]$.
We choose $\lambda_1 = \lambda_2 = 0.3$ and $\lambda_3 = \lambda_4 = 0.2$.

The functions $r_{\text{gap}}$, $r_{\text{eff}}$, and $r_{\text{conf}}$ are defined as follows.
Assuming that $h$ was successful, let us denote by $x_{new}$ and $x_{old}$ the new and old solution, respectively.
Furthermore, let $x_{LP}$ be the solution of the current linear relaxation.
Then, we measure the quality of $x_{new}$ relative to the current solving stage with
\[
r_{\text{gap}}(h,t) = 
\begin{cases}
    0, & \text{if } h \text{ did not find an incumbent at iteration } t \\
    1, & \text{if } h \text{ found the first incumbent at iteration } t \\
    \frac{c^T x_{old} - c^T x_{new}}{c^T x_{old} - c^T x_{LP}}, & \text{otherwise.}
\end{cases}
\]
To define $r_{\text{eff}}$, let $n_h^t$ be the number of nodes used by $h$, and 
$n_{max}$ an upper bound on the maximal number of nodes used.
For LNS, $n_h^t$ refers to the number of nodes solved in the sub-MIP; for diving, it refers to the number of nodes visited during the partial search.
Finally, we define
$
r_{\text{eff}}(h,t) = 1 - \frac{n_h^t}{n_{max}}
$
and
$
r_{\text{conf}}(h,t) = \frac{v_h^t}{v_{max}}
$
where $v_h^t$ is the number of conflict constraints $h$ found and $v_{max}$ the maximal number of conflict constraints found by any heuristic in the past.
The reward function $r$ is an extension of the reward used in \cite{hendel18}, which only uses $r_{\text{gap}}$ and variants of $r_{\text{sol}}$ and $r_{\text{eff}}$.

\subsection{Choosing a Bandit Algorithm}
\label{sec:bandit}

As mentioned before, to solve the multi-armed bandit problem successfully, we need to balance exploitation and exploration carefully.
In our case, this raises the following question:
Should we prioritize heuristics that have not been executed (that often) or heuristics that have performed well in the past?
Our experimental results suggests that for primal heuristics, exploitation is the better choice.
Typically, a heuristic that performs bad at the beginning, will also be rather unsuccessful later on, since it only gets harder to find improving solutions.

That is why we propose to use a variant of the $\epsilon$-greedy bandit algorithm. The $\epsilon$-greedy, or follow-the-leader, algorithm pursues a simple strategy:
Given an $\epsilon \in [0,1]$, the best action seen so far is chosen with probability $1 - \epsilon$.
Otherwise, an action is randomly selected following a uniform distribution.
To characterize the best action at iteration $t$, we associated a weight $w(h,t)$ with every heuristic~$h$.
The weights $w$ are equal to the average reward of $h$ observed so far, that is,
$w(h,t) = \sum_{\tilde{t} \in \mathcal{T}_h^t} r(h, \tilde{t}) / |\mathcal{T}_h^t|$
with $\mathcal{T}_h^t \subseteq [t]$ being the subset of calls at which $h$ was selected up to time $t$. 

In the modified $\epsilon$-greedy algorithm we consider, instead of selecting a heuristic uniformly at random, we draw it following the distribution imposed by the weights $w$.
This variant allows for more exploitation; it is described in Algorithm~\ref{alg:greedy}.
In our experiments, our online scheduling approach performed best with this bandit algorithm.
We use $\epsilon = 0.7$.

To put more focus on heuristics that performed well in the \emph{recent} past, we also tried another modification:
Instead of looking at the average reward as $w$, we examined using an aggregation of the observed rewards where older observations contribute exponentially less.
This performed considerably worse, suggesting that it is preferable to consider all past behavior to make future decisions.

\begin{algorithm}[t]
	\caption{Modified $\epsilon$-greedy bandit algorithm}
	\label{alg:greedy}
	\begin{algorithmic}
		\STATE {\bfseries Input:} Set of heuristics $\mathcal{H}$, reward function $r$, probability $\epsilon \in [0,1]$

        \STATE $w(h, 0) \gets \frac{1}{|\mathcal{H}|}$
        \STATE $t \gets 0$

		\WHILE{not stopped}
		    \STATE $t \gets t + 1$
		    \STATE $\epsilon_t \gets \epsilon \cdot \sqrt{\frac{|\mathcal{H}|}{t}}$
		    \STATE Draw $\rho_t \sim \mathbb{U}([0,1])$
		    \IF{$\rho_t > \epsilon_t$}
		        \STATE $h_t \gets \underset{h \in \mathcal{H}}{\text{argmax}}$ $w(h,t-1)$
		    \ELSE
		        \STATE Draw $h_t \sim w(\cdot, t-1)$
		    \ENDIF
		
		\STATE Observe reward $r(h_t, t)$
		
		\IF{$h_t$ was selected for the first time}
		    \STATE $w(h_t, t) \gets r(h_t, t)$
		\ELSE
		    \STATE $w(h_t, t) \gets$ update average weight with $r(h_t, t)$
		\ENDIF
		
		\ENDWHILE
		
	\end{algorithmic}
\end{algorithm}

\section{Computational Results}
\label{sec:results}

To study the performance of our approach, we used the state-of-the-art open-source MIP solver SCIP 8.0 with SoPlex 6.0 \cite{SCIP}.
We ran all experiments on a Linux cluster of Intel Xeon CPU E5-2630 v3 2.40GHz with 64GB RAM.
The time limit was set two hours and the test set consists of the benchmark instances of the MIPLIB 2017 \cite{miplib}.
Since our framework aims to improve primal performance, we removed all infeasible instances and problems with a zero objective function.
This leaves us with 226 instances.
To filter out the effects of performance variability \cite{lodi13}, all experiments are run with four random seeds.

We compare two settings:
\textsc{default} refers to the default settings of SCIP and \textsc{scheduler} refers to the proposed online scheduling framework.
In the latter, we deactivated all LNS and diving heuristics that are controlled by the scheduler, as well as the two adaptive heuristics presented in \cite{hendel18, hendel182}.
The scheduler is called at every node, right after cheaper heuristics like rounding.

\noindent
Table~\ref{tab:bbresults} shows that \textsc{scheduler} consistently reduces the size of the B\&B tree by $4-6\%$.
Unfortunately, this improvement does not directly translate into improving 
solving time.
However, we perform the better the harder the instances get:
On instances taking at least $1000$ seconds to solve, the scheduling framework outperforms \textsc{default} by about $4\%$.
Even though \textsc{scheduler} solves 13 instances that cannot be solved by \textsc{default}, it fails to solve 21 instances solved by \textsc{default}.
One reason for this behavior could be the fact that for harder instances, the scheduler tends to spend less time then default in the heuristics controlled by it:
On $[1000,7200]$, \textsc{scheduler} reduces time spent in heuristics by over $30\%$.

The heuristics' behavior shows that our scheduling framework succeeds in detecting more successful heuristics:
The scheduler finds $88\%$ more incumbents while increasing the probability of a heuristic finding a new solution by $57\%$.
On average, the heuristics controlled by \textsc{scheduler} find 3.80 solutions per instance as opposed to 2.01; with a success probability of $4.49\%$ instead of $2.86\%$.

To conclude, the computational results show that our scheduling framework can improve the performance of a solver.
However, for easier instances, it seems that there is not much potential for 
improvement by using an online learning approach; there we compete against the 
default parameters that have been tuned on the test set over a long period of 
time.
This could be attributed to a lack of observations:
When an instance is solved fast, a learning approach might not have enough time to gather meaningful information about the heuristics' behavior.  
Furthermore, the results also suggest that our approach might be too conservative for harder instances, since it reduces the time spent in heuristics considerably.

Hence, we believe that there is further room for improvement, also since we have not spent a large amount of effort on tuning any hyperparameters of our method in order to obtain the current results.
As next steps, we need to combine the good performance of the static heuristic handling with our online scheduling approach and better detect when to apply heuristics more aggressively and when to rely on well-working default parameters.

\begin{table}[t]
\caption{Summary of results for B\&B experiments. Rows labeled $[t,7200]$ consist of instances solved with at least one settings taking at least $t$ seconds. \textit{heurtime} refers to time spent in heuristics controlled by the scheduler, relative to \textsc{default}. Shifted geometric means are used.} 
\label{tab:bbresults}
\renewcommand{\arraystretch}{1.35}
\scriptsize
\centering
\resizebox{\columnwidth}{!}{%
\begin{tabular}{lccccccccccccc}
\hline
\hline
& & & \multicolumn{3}{c}{\textsc{default}} && \multicolumn{3}{c}{\textsc{scheduler}} && \multicolumn{3}{c}{relative} \\
\cline{4-6} \cline{8-10} \cline{12-14}
subset                & instances &&  solved &        time &        nodes &&     solved &       time &         nodes &&        time &         nodes &   heurtime\\
\hline
all                   &       892 &&     472 &     1157.44 &         4238 &&        464 &     1189.03 &         4056 &&        1.03 &          0.95 &       0.94\\
\hline
$[0,7200]$            &       485 &&     472 &      249.02 &         2522 &&        464 &      261.70 &         2426 &&        1.05 &          0.96 &       1.20\\
$[1,7200]$            &       481 &&     468 &      259.89 &         2593 &&        460 &      273.23 &         2494 &&        1.05 &          0.96 &       1.20\\
$[10,7200]$           &       441 &&     428 &      373.99 &         3439 &&        420 &      394.12 &         3298 &&        1.05 &          0.96 &       1.21\\
$[100,7200]$          &       330 &&     317 &      839.49 &         8231 &&        309 &      862.90 &         7759 &&        1.03 &          0.94 &       1.05\\
$[1000,7200]$         &       175 &&     162 &     2312.56 &        20769 &&        154 &     2217.32 &        19627 &&        0.96 &          0.95 &       0.68\\
\hline
all-optimal           &       451 &&     451 &      199.60 &         2294 &&        451 &      209.61 &         2168 &&        1.05 &          0.95 &       1.25\\
\hline
\hline
\end{tabular}
}
\end{table}

%
%
%
 \bibliographystyle{splncs04}
 \bibliography{biblio}

\begin{thebibliography}{10}
\providecommand{\url}[1]{\texttt{#1}}
\providecommand{\urlprefix}{URL }
\providecommand{\doi}[1]{https://doi.org/#1}

\bibitem{achterberg07}
Achterberg, T.: Conflict analysis in mixed integer programming. Discrete
  Optimization  \textbf{4}(1),  4--20 (2007)

\bibitem{balcan2018learning}
Balcan, M.F., Dick, T., Sandholm, T., Vitercik, E.: Learning to branch. In:
  International conference on machine learning. pp. 344--353. PMLR (2018)

\bibitem{baltean19}
Baltean-Lugojan, R., Bonami, P., Misener, R., Tramontani, A.: Scoring positive
  semidefinite cutting planes for quadratic optimization via trained neural
  networks. preprint: \url{https://optimization-online.org/?p=17362}  (2019)

\bibitem{berthold132}
Berthold, T.: Measuring the impact of primal heuristics. Operations Research
  Letters  \textbf{41}(6),  611--614 (2013)

\bibitem{berthold13}
Berthold, T.: Primal {MINLP} heuristics in a nutshell. In: International
  Conference on Operations Research (2013)

\bibitem{berthold18}
Berthold, T.: A computational study of primal heuristics inside an {MI(NL)P}
  solver. Journal of Global Optimization  \textbf{70},  189--206 (01 2018)

\bibitem{berthold17}
Berthold, T., Hendel, G., Koch, T.: From feasibility to improvement to proof:
  three phases of solving mixed-integer programs. Optimization Methods and
  Software  \textbf{33},  1--19 (11 2017)

\bibitem{SCIP}
Bestuzheva, K., Besan{\c{c}}on, M., Chen, W.K., Chmiela, A., Donkiewicz, T.,
  van Doornmalen, J., Eifler, L., Gaul, O., Gamrath, G., Gleixner, A.,
  Gottwald, L., Graczyk, C., Halbig, K., Hoen, A., Hojny, C., van~der Hulst,
  R., Koch, T., L{\"u}bbecke, M., Maher, S.J., Matter, F., M{\"u}hmer, E.,
  M{\"u}ller, B., Pfetsch, M.E., Rehfeldt, D., Schlein, S., Schl{\"o}sser, F.,
  Serrano, F., Shinano, Y., Sofranac, B., Turner, M., Vigerske, S.,
  Wegscheider, F., Wellner, P., Weninger, D., Witzig, J.: {The SCIP
  Optimization Suite 8.0}. ZIB-Report 21-41, Zuse Institute Berlin (December
  2021), \url{http://nbn-resolving.de/urn:nbn:de:0297-zib-85309}

\bibitem{bubeck12}
Bubeck, S., Nicolò, C.B.: Regret Analysis of Stochastic and Nonstochastic
  Multi-armed Bandit Problems. Foundations and Trends In Machine Learning
  (2012)

\bibitem{chmiela21}
Chmiela, A., Khalil, E., Gleixner, A., Lodi, A., Pokutta, S.: Learning to
  schedule heuristics in branch and bound. In: Advances in Neural Information
  Processing Systems. vol.~34 (2021)

\bibitem{dambrosio20}
D'Ambrosio, C., Antonio, F., Iommazzo, G., Liberti, L.: A learning-based
  mathematical programming formulation for the automatic configuration of
  optimization solvers. Machine Learning, Optimization, and Data Science pp.
  700--712 (01 2020)

\bibitem{etheve20}
Etheve, M., Alès, Z., Bissuel, C., Juan, O., Kedad-Sidhoum, S.: Reinforcement
  learning for variable selection in a branch and bound algorithm.
  arXiv:2005.10026  (2020)

\bibitem{miplib}
Gleixner, A., Hendel, G., Gamrath, G., Achterberg, T., Bastubbe, M., Berthold,
  T., Christophel, P.M., Jarck, K., Koch, T., Linderoth, J., L\"ubbecke, M.,
  Mittelmann, H.D., Ozyurt, D., Ralphs, T.K., Salvagnin, D., Shinano, Y.:
  {MIPLIB 2017: Data-Driven Compilation of the 6th Mixed-Integer Programming
  Library}. Mathematical Programming Computation  (2021)

\bibitem{he14}
He, H., III, H.D., Eisner, J.M.: Learning to search in branch and bound
  algorithms. In: Advances in Neural Information Processing Systems. vol.~27,
  pp. 3293--3301 (2014)

\bibitem{hendel18}
Hendel, G.: Adaptive large neighborhood search for mixed integer programming.
  Mathematical Programming Computation  \textbf{14}(2),  185--221 (2022)

\bibitem{hendel182}
Hendel, G., Miltenberger, M., Witzig, J.: Adaptive algorithmic behavior for
  solving mixed integer programs using bandit algorithms. In: International
  Conference on Operations Research (2018)

\bibitem{chen22}
Huang, L., Chen, X., Huo, W., Wang, J., Zhang, F., Bai, B., Shi, L.: Improving
  primal heuristics for mixed integer programming problems based on problem
  reduction: A learning-based approach. arXiv:2209.13217  (2022)

\bibitem{huang22}
Huang, Z., Wang, K., Liu, F., Zhen, H.L., Zhang, W., Yuan, M., Hao, J., Yu, Y.,
  Wang, J.: Learning to select cuts for efficient mixed-integer programming.
  Pattern Recognition  \textbf{123} (2022)

\bibitem{hutter09}
Hutter, F., Hoos, H., Leyton-Brown, K., St{\"u}tzle, T.: Paramils: An automatic
  algorithm configuration framework. Journal of Artificial Intelligence
  Research (JAIR)  \textbf{36},  267--306 (10 2009)

\bibitem{hutter11}
Hutter, F., Hoos, H.H., Leyton-Brown, K.: Sequential model-based optimization
  for general algorithm configuration. In: Learning and Intelligent
  Optimization. pp. 507--523 (2011)

\bibitem{khalil16}
Khalil, E.B., Bodic, P.L., Song, L., Nemhauser, G., Dilkina, B.: Learning to
  branch in mixed integer programming. In: Proceedings of the 30th AAAI
  Conference on Artificial Intelligence (2016)

\bibitem{kruber17}
Kruber, M., L{\"u}bbecke, M., Parmentier, A.: Learning when to use a
  decomposition. In: International Conference on AI and OR Techniques in
  Constraint Programming for Combinatorial Optimization Problems. pp. 202--210
  (2017)

\bibitem{land60}
Land, A.H., Doig, A.G.: An automatic method of solving discrete programming
  problems. Econometrica  \textbf{28}(3),  497--520 (1960)

\bibitem{lattimore20}
Lattimore, T., Szepesvári, C.: Bandit Algorithms. Cambridge University Press
  (2020)

\bibitem{lodi13}
Lodi, A., Tramontani, A.: Performance variability in mixed-integer programming.
  Tutorials in Operations Research, Vol. 10 pp. 1--12 (09 2013)

\bibitem{lodi17}
Lodi, A., Zarpellon, G.: On learning and branching: a survey. TOP: An Official
  Journal of the Spanish Society of Statistics and Operations Research
  \textbf{25},  207--236 (2017)

\bibitem{nair20}
Nair, V., Bartunov, S., Gimeno, F., von Glehn, I., Lichocki, P., Lobov, I.,
  O'Donoghue, B., Sonnerat, N., Tjandraatmadja, C., Wang, P., Addanki, R.,
  Hapuarachchi, T., Keck, T., Keeling, J., Kohli, P., Ktena, I., Li, Y.,
  Vinyals, O., Zwols, Y.: Solving mixed integer programs using neural networks.
  arXiv preprint: 2012.13349  (2020)

\bibitem{paulus22}
Paulus, M.B., Zarpellon, G., Krause, A., Charlin, L., Maddison, C.: Learning to
  cut by looking ahead: Cutting plane selection via imitation learning. In:
  Proceedings of the 39th International Conference on Machine Learning.
  vol.~162, pp. 17584--17600 (2022)

\bibitem{scavuzzo22}
Scavuzzo, L., Chen, F.Y., Chételat, D., Gasse, M., Lodi, A., Yorke-Smith, N.,
  Aardal, K.: Learning to branch with tree mdps. arXiv:2205.11107  (2022)

\bibitem{tang20}
Tang, Y., Agrawal, S., Faenza, Y.: Reinforcement learning for integer
  programming: Learning to cut. In: Proceedings of the 37th International
  Conference on Machine Learning. vol.~119, pp. 9367--9376 (2020)

\bibitem{turner22}
Turner, M., Koch, T., Serrano, F., Winkler, M.: Adaptive cut selection in
  mixed-integer linear programming. arXiv:2202.10962  (2022)

\bibitem{yilmaz21}
Yilmaz, K., Yorke-Smith, N.: A study of learning search approximation in mixed
  integer branch and bound: Node selection in scip. AI  \textbf{2},  150--178
  (04 2021)

\end{thebibliography}
\end{document}